\documentclass[12pt,a4paper]{amsart}
\usepackage{amssymb}
\usepackage{mathrsfs}
\headsep=1cm \numberwithin{equation}{section}
\newtheorem{thm}{Theorem}[section]
\newtheorem{cor}[thm]{Corollary}
\newtheorem{lem}[thm]{Lemma}
\newtheorem{ex}[thm]{Example}
\newtheorem{prop}[thm]{Proposition}
\theoremstyle{definition}

\theoremstyle{remark}
\newtheorem{rem}[thm]{Remark}
\numberwithin{equation}{section}

\newcommand\Supp{\operatorname{Supp}}
\newcommand\Ass{\operatorname{Ass}}
\newcommand\mAss{\operatorname{mAss}}

\newcommand\Rad{\operatorname{Rad}}
\newcommand\Hom{\operatorname{Hom}}
\newcommand\Ext{\operatorname{Ext}}
\newcommand\Tor{\operatorname{Tor}}

\newcommand\height{\operatorname{height}}

\newcommand\ara{\operatorname{ara}}
\DeclareMathOperator{\Ker}{ker} \DeclareMathOperator{\Coker}{coker}
\DeclareMathOperator{\lc}{H} \DeclareMathOperator{\Image}{im}
\newcommand\cd{\operatorname{cd}}
\newcommand\m{\operatorname{\frak m}}
\newcommand\q{\operatorname{\frak q}}
\newcommand\p{\operatorname{\frak p}}

\newcommand{\lo}{\longrightarrow}

\begin{document}\title[Some homological properties of ideals]
{Some homological properties of ideals with cohomological dimension one}
\author[G. Pirmohammadi, K. Ahmadi Amoli, and K. Bahmanpour]{G. Pirmohammadi, K. Ahmadi Amoli, and K. Bahmanpour$^*$ }

\address{Gholamreza Pirmohammadi;
 Payame Noor University, Po Box 19395-3697, Tehran, Iran.} \email{\it pirmohammadi.reza@gmail.com}
\address{Khadijeh Ahmadi Amoli;
 Payame Noor University, Po Box 19395-3697, Tehran, Iran.}  \email {\it khahmadi@pnu.ac.ir}
\address{Kamal Bahmanpour; Faculty of Sciences, Department of Mathematics, University of Mohaghegh Ardabili,
56199-11367, Ardabil, Iran;
and School of Mathematics, Institute for Research in Fundamental Sciences (IPM), P.O. Box. 19395-5746, Tehran, Iran.} \email{\it bahmanpour.k@gmail.com}

\thanks{ $^*$Corresponding author: e-mail:{\it bahmanpour.k@gmail.com} (Kamal Bahmanpour)}%

\begin{abstract}
Let $R$ denote a commutative Noetherian ring and let $I$ be an ideal of $R$ such that $H^i_I(R)=0$,
 for all integers $i\geq 2$.  In this paper we shall prove some results concerning the homological properties of $I$.
\end{abstract}
\maketitle
\section{Introduction}
Let $R$ denote a commutative Noetherian ring and  let $I$ be an
ideal of $R$. In \cite{Ha}, Hartshorne defined an $R$-module $L$ to be
$I$-{\it cofinite}, if $\Supp L\subseteq
V(I)$ and ${\rm Ext}^{i}_{R}(R/I, L)$ is a finitely generated module
for all $i$. He posed the following question:\\

{\it Whether the category $\mathscr{C}(R, I)_{cof}$
of $I$-cofinite modules is an Abelian subcategory of the category of all $R$-modules?
That is, if $f: M\longrightarrow N$ is an $R$-homomorphism of
$I$-cofinite modules, are $\ker f$ and ${\rm coker} f$ $I$-cofinite?}\\

Hartshorne gave a counterexample to show that this question has not an
affirmative answer in general, (see \cite[\S 3]{Ha}). On the positive side, Hartshorne proved that if $I$ is a prime
ideal of dimension one in a complete regular local ring $R$, then the answer to his question
is yes. On the other hand, in \cite{DM}, Delfino and Marley extended this result to arbitrary
complete local rings. Kawasaki in \cite{Ka2} generalized the Delfino and Marley's result
for an arbitrary ideal $I$ of dimension one in a local ring $R$. Finally, Melkersson
in \cite{Me2} generalized the Kawasaki's result for all ideals of dimension one of any
 arbitrary Noetherian ring $R$. More recently, in \cite{BNS1} as a generalization of
 Melkersson's result it is shown that for any ideals $I$ in a commutative Noetherian ring $R$, the category of all $I$-cofinite
$R$-modules $M$ with $\dim M\leq 1$ is an Abelian subcategory of the category of all $R$-modules.\\

 Recall that, for an $R$-module $M$, the {\it cohomological dimension of $M$ with
respect to} $I$, denoted by $\cd (I, M)$, is the smallest integer $n\geq 0$ such that $H^i_I(M)=0$ for all $i>n$.\\

The cohomological dimension have been studied by
several authors; see, for example, Faltings \cite{Fa}, Hartshorne
\cite{Ha1}, Huneke-Lyubeznik \cite{HL}, Divaani-Aazar et al \cite{DNT}, Hellus \cite{He}, Hellus-St\"uckrad \cite{HS},
Mehrvarz et al \cite{MBN},  and Ghasemi et al \cite{GBA}.\\

Recall that, for any proper ideal $I$ of $R$, the
\emph{arithmetic rank} of $I$, denoted by $\ara(I)$, is the least number of
elements of $I$ required to generate an ideal which has the same radical as $I$.\\

Now, let $I$ be an ideal of an arbitrary Noetherian ring $R$. Kawasaki in \cite[Theorem 2.1]{Ka1}
 proved that if $\ara(I)=1$ then the category $\mathscr{C}(R, I)_{cof}$
of $I$-cofinite modules is an Abelian subcategory of the category of all $R$-modules.\\

It is well known that, for any proper ideal $I$ of a Noetherian ring $R$ we have $\cd (I, R)\leq\ara(I)$.
 (See \cite[Corollary 3.3.3]{BS}). In particular, for any ideal $I$ with $\ara(I)=1$ we have $\cd (I, R)\leq 1$.
  So, as a generalization of Kawasaki's interesting result, it is more natural that we ask the following question:\\

Question 1: {\it Let $R$ be a Noetherian ring and $I$ be an ideal of $R$ with $\cd (I, R)\leq 1$.
Whether the category $\mathscr{C}(R, I)_{cof}$
of $I$-cofinite modules is an Abelian subcategory of the category of all $R$-modules?}\\

In section 2 of this paper we present an affirmative answer to the Question 1, whenever $R$ is a local Noetherian ring.
 In section 3 we present some equivalent conditions for the exactness of ideal transform functors. Using these results,
  for any ideal $I$ generated by two elements,  we provide some necessary and sufficient conditions for the non-vanishing
  of the local cohomology module $H^2_I(R)$. Finally, in section 4 we prepare some vanishing conditions for the Bass and
  Betti numbers of some special local cohomology modules. Also, we give a formula for the cohomological dimension of some
  special ideals in Noetherian domains.\\\\

Throughout this paper, $R$ will always be a commutative Noetherian
ring and $I$ will be an ideal of $R$. Also, for an $R$-module $M$, $\Gamma _I(M)$ denotes the submodule of $M$ consisting
 of all elements annihilated by
some power of $I$, i.e., $\bigcup_{n=1}^{\infty}(0:_MI^n)$. For every $R$-module $L$, we denote by $\mAss_R L$ the set of
 minimal elements of $\Ass_R L$ with respect to inclusion. Also, for any ideal $\frak a$ of $R$, we denote
$\{\frak p \in {\rm Spec}\,R:\, \frak p\supseteq \frak a \}$ by
$V(\frak a)$. Finally, for any ideal $\frak{b}$ of $R$, {\it the
radical of} $\frak{b}$, denoted by $\Rad(\frak{b})$, is defined to
be the set $\{x\in R \,: \, x^n \in \frak{b}$ for some $n \in
\mathbb{N}\}$.  For any unexplained notation and terminology we refer
the reader to \cite{BS} and \cite{Mat}.

\section{A category of modules which is Abelian}

The following well known lemma plays a key role in the proof of Proposition 2.2.

\begin{lem}
\label{2.1}
Let $R$ be a Noetherian ring, $I$ be an ideal of $R$ and $T$ be an $R$-module. Then  the following statements are equivalent:
\begin{itemize}
 \item[(i)] The $R$-module $H^i_I(T)$ is Artinian, for $i\geq 0$.
 \item[(ii)] The $R$-module $\Ext^i_R(R/I,T)$ is Artinian, for $i\geq 0$.
\item[(iii)] The $R$-module $\Ext^i_R(R/I,T)$ is Artinian, for $0 \leq i \leq \cd(I,R)$.
\end{itemize}
\end{lem}
\proof See \cite[Theorem 5.5]{Me} and \cite[Theorem 2.9]{AM}.\qed\\

\begin{prop}
\label{2.2}
 Let $(R,\m)$ be a Noetherian local ring, $I$ be an ideal of $R$ and $M$ be an $R$-module. Then the following conditions
are equivalent:
\begin{itemize}
 \item[(i)] The $R$-modules $\Tor^R_{i}(R/I,M)$ are finitely generated for all $i\geq 0$.
 \item[(ii)] The $R$-modules $\Tor^R_{i}(R/I,M)$ are finitely generated for all $0\leq i\leq \cd(I,R)$.
\end{itemize}
\end{prop}
\proof
(i)$\Rightarrow$(ii) Is clear.\\

(ii)$\Rightarrow$(i) Let $T:=D(M)$, where $D(-):= {\rm Hom}_{R}(-,E)$
denotes the Matlis dual functor and $E := E_R(R/\frak{m})$ is the
injective hull of the residue field $R/\frak{m}$. Then by the adjointness we have $\Ext^i_R(R/I,T)\simeq D(\Tor^R_i(R/I,M))$,
 for all $0\leq i\leq \cd(I,R)$. In particular, the $R$-modules $\Ext^i_R(R/I,T)$ are Artinian for all $0 \leq i \leq \cd(I,R)$.
  So, it follows from Lemma 2.1 and adjointness, that the $R$-modules $\Ext^i_R(R/I,T)\simeq D(\Tor^R_i(R/I,M))$ are
  Artinian for all $i\geq 0$. Now, it follows from \cite[Lemma 1.15(a)]{KLW} that, the $R$-modules $\Tor^R_i(R/I,M)$
  are finitely generated for all $i\geq 0$,
as required.\qed\\

The following easy consequence of Proposition 2.2 plays a key role in the proof of Theorem 2.4.

\begin{cor}
\label{2.3}
Let $(R,\m)$ be a Noetherian local ring, $I$ be an ideal of $R$ with $\cd(I,R)=n$
and let $M$ be an $R$-module with $\Supp M\subseteq V(I)$. Then  the $R$-module $M$ is $I$-cofinite,
 if and only if, the $R$-modules $\Tor^R_i(R/I,M)$  are finitely generated for all $0 \leq i \leq n$.
\end{cor}
\proof The assertion follows from Proposition 2.2 and \cite[Theorem 2.1]{Me}.\qed\\

The following result gives an affirmative answer to the Question 1, for the local case.

\begin{thm}
 \label{2.4}
Let $I$ be an ideal of a Noetherian local ring $(R,\m)$ such that $\cd(I,R)\leq 1$. Let
$\mathscr{C}(R, I)_{cof}$ denote the category of $I$-cofinite
$R$-modules. Then
$\mathscr{C}(R, I)_{cof}$ is an Abelian subcategory of the category of all $R$-modules.
\end{thm}
\proof Let $M,N\in \mathscr{C}(R, I)_{cof}$ and let $f:M\longrightarrow N$ be
an $R$-homomorphism. It is enough to show that the $R$-modules
$\ker f$ and ${\rm coker} f$ are $I$-cofinite.

To this end, the exact sequence $$0\longrightarrow \ker f \longrightarrow M \longrightarrow
{\rm im}f \longrightarrow 0,$$ induces an exact sequence $$\Tor^R_0(R/I,M)\rightarrow
\Tor^R_0(R/I,{\rm im}f)\rightarrow 0,$$which using \cite[Theorem 2.1]{Me}, implies that
 the $R$-module $\Tor^R_0(R/I,{\rm im}f)$ is finitely generated. Now, the exact sequence
$$0\longrightarrow {\rm im}f \longrightarrow N \longrightarrow
{\rm coker}f \longrightarrow 0.$$
induces an  exact sequence

$$\Tor^R_1(R/I,N)\rightarrow \Tor^R_1(R/I,{\rm coker}f)\rightarrow \Tor^R_0(R/I,{\rm im}f)$$$$
\rightarrow  \Tor^R_0(R/I,N)\rightarrow  \Tor^R_0(R/I,{\rm coker}f)\rightarrow 0.$$ By \cite[Theorem 2.1]{Me}
 the modules $\Tor^R_1(R/I,N)$ and $\Tor^R_0(R/I,N)$ are finitely generated $R$-modules, which implies that the
$R$-modules $\Tor^R_0(R/I,{\rm coker}f)$ and $\Tor^R_1(R/I,{\rm coker}f)$ are finitely generated. Therefore, it follows
from Corollary 2.3, that the
$R$-module ${\rm coker}f$ is $I$-cofinite. Now, the assertion follows from the following exact sequences
$$0\longrightarrow {\rm im}f \longrightarrow N \longrightarrow
{\rm coker}f \longrightarrow 0,$$ and $$0\longrightarrow \ker f \longrightarrow M \longrightarrow
{\rm im} f \longrightarrow 0.$$
\qed \\

\begin{cor}
 \label{2.5}
Let $I$ be an ideal of a Noetherian local ring $(R,\m)$ such that $\cd(I,R)\leq1$.
Let $\mathscr{C}(R, I)_{cof}$ denote the category of $I$-cofinite modules over $R$.
Let $$X^\bullet:
\cdots\longrightarrow X^i \stackrel{f^i} \longrightarrow X^{i+1}
\stackrel{f^{i+1}} \longrightarrow X^{i+2}\longrightarrow \cdots,$$ be a
complex such that  $X^i\in\mathscr{C}(R, I)_{cof}$ for all $i\in\Bbb{Z}$.
Then for each $i\in\Bbb{Z}$ the $i^{{\rm th}}$ cohomology module $H^i(X^\bullet)$ is in
$\mathscr{C}(R, I)_{cof}$.
\end{cor}
\proof The assertion follows from  Theorem 2.4.\qed\\

\begin{cor}
\label{2.6}
Let $(R,\m)$ be a Noetherian local ring, $I$ be an ideal of $R$ with $\cd(I,R)\leq1$ and
let $M$ be an $I$-cofinite $R$-module. Then,  the $R$-modules
${\rm Tor}_i^R(N,M)$ and ${\rm Ext}^i_R(N,M)$ are $I$-cofinite, for all finitely generated $R$-modules $N$ and all
integers $i\geq0$.
\end{cor}
\proof Since $N$ is finitely generated it follows that, $N$ has a
free resolution with finitely generated free $R$-modules. Now the
assertion follows using Corollary 2.5 and computing the $R$-modules ${\rm
Tor}_i^R(N,M)$ and ${\rm Ext}^i_R(N,M)$, using this
free resolution. \qed\\

\section{Vanishing of the extension and torsion functors}

In this section we present some equivalent conditions for the exactness of ideal transform functors.
Using these results, for any  ideal $I$ generated by two elements,  we provide some
necessary and sufficient conditions for the non-vanishing of the local cohomology module $H^2_I(R)$.\\

The following lemma is needed in the proof of Proposition 3.2.

\begin{lem}
\label{3.1}
 Let $R$ be a Noetherian ring and $I$ be an
ideal of $R$. Let $M$ be an $R$-module such that $\Tor^R_{i}(R/I,M)=0$, for all integers $i\geq 0$. Then
  $\Hom_{R}(R/I,M)=0$.
\end{lem}
\proof
Let $I=(x_1,\dots,x_n)$ and let
$$K_\bullet(\underline{x}, M): 0\lo M \overset{f_0}\lo  \bigoplus_{k=1}^{C_n^1} M
\overset{f_1}\lo \bigoplus_{k=1}^{C_n^2} M \lo \dots \lo
\bigoplus_{k=1}^{C_n^{n-1}} M\overset{f_{n-1}} \lo M\lo 0,$$
be the Koszul complex of $M$ with respect to $\underline{x}=x_1,\dots,x_n$.\\

We prove that $\lc_{i}(\underline{x}; M)=0$ for all $0\leq i \leq n$. By the definition we have
$$\lc_{0}(\underline{x}; M)= \Coker f_{n-1}= M/IM\simeq R/I\otimes_{R}M=0.$$
So, we have $\Image f_{n-1}=M$.
Therefore, using the hypothesis it follows  that
$$\Tor^R_i(R/I,\Image f_{n-1})=0,\,\,{\rm for}\,\,{\rm each}\,\, i\geq 0.$$ Hence,
the exact sequence

\begin{center}
$0\lo \Ker{ f_{n-1}}\lo \bigoplus_{k=1}^{C_n^{n-1}} M\lo \Image f_{n-1} \lo 0,$
\end{center}
implies that $\Tor^R_i(R/I,\Ker f_{n-1})=0$, for each $i\geq 0$.
The exact sequence
\begin{center}
$0\lo \Image{ f_{n-2}}\lo \Ker f_{n-1}\lo \lc_1(\underline{x}; M)\lo 0,
\,\,\,\,\,\,~(3.1.1)$
\end{center}
induces the exact sequence
\begin{center}
$R/I\otimes_{R}\Ker f_{n-1}  \lo R/I\otimes_{R} \lc_1(\underline{x};M)\lo
0.\,\,\,\,\,\,~(3.1.2)$
\end{center}
Now as $R/I\otimes_{R}\Ker f_{n-1} =0$, it follows from $(3.1.2)$ that
 $$R/I\otimes_{R} \lc_1(\underline{x};M)=0.$$ By
the definition of the Koszul complex we have $I\lc_1(\underline{x};M)=0$. Therefore, we have $\lc_1(\underline{x};M)\simeq R/I\otimes_{R} \lc_1(\underline{x};M)=0$. Now it follows from the exact sequence $(3.1.1)$ that $\Tor^R_i(R/I,\Image f_{n-2})=0$, for each $i\geq 0$.
Moreover, the exact sequence

\begin{center}
$0\lo \Ker f_{n-2}  \lo \bigoplus_{k=1}^{C_n^{n-2}} M \lo \Image{ f_{n-2}}\lo 0,$
\end{center}
implies that
$\Tor^R_i(R/I,\Ker f_{n-2})=0$, for each $i\geq 0$. Proceeding in the same way we can see
$\lc_i(\underline{x};M)=0$, for all $0\leq i \leq n$.\\

Now, since
 $\lc_i(\underline{x};M)\simeq \lc^{n-i}(\underline{x};M)$ for all $0\leq i \leq n$, it follows
 that $$\Hom_R(R/I,M)\simeq(0:_MI)\simeq \lc^{0}(\underline{x};M)\simeq \lc_n(\underline{x};M)=0.$$\qed\\

The following proposition plays a key role in the proof of our main results.

\begin{prop}
\label{3.2}
  Let $R$ be a Noetherian ring, $I$ be an
ideal of $R$ and $M$ be an $R$-module. Then the following conditions
are equivalent:
\begin{itemize}
 \item[(i)]  $\Tor^R_{n}(R/I,M)=0$, for all integers $n\geq 0$.
 \item[(ii)] $\Ext^n_{R}(R/I,M)=0$, for all integers $n\geq 0$.
 \end{itemize}
\end{prop}
\proof
(i)$\Rightarrow$(ii) We argue by induction on $n$. For $n=0$, the assertion holds by Lemma 3.1.
We therefore assume, inductively, that $n>0$ and the result has been proved for smaller values of $n$.
 Then there is an exact sequence $$0  \longrightarrow M  \longrightarrow E_R(M)  \longrightarrow E_R(M)/M  \longrightarrow0.\,\,\,\,\,\,~(3.2.1)  $$
Since, by the hypothesis we have $(0:_MI)=0$ it follows that $(0:_{E_R(M)}I)=0$ and hence $\Gamma_I(E_R(M))=0$.
By \cite[Theorem 18.5]{Mat} the $R$-module $E_R(M)$ is isomorph with a direct sum of a family of indecomposable injective $R$-modules.
 Let $\p$ be a prime ideal of $R$ such that $E_R(R/\p)$ is a direct summand of  $E_R(M)$. Then we have $\Gamma_I(E_R(R/\p))=0$.
 Therefore, form the fact that $\Ass_R E_R(R/\p)=\{\p\}$ it follows that $I\not\subseteq \p$.
 So there exists an element $a\in I$ such that $a\not\in\p$. Then by \cite[Theorem 18.4(iii)]{Mat},
 multiplication by $a$ is an automorphism on $E_R(R/\p)$.  Therefore, multiplication by $a$ is an automorphism on
  $\Tor^R_i(R/I,E_R(R/\p))$, for all $i\geq 0$. But, since $a\in I$ it follows that, multiplication by $a$ on $\Tor^R_i(R/I,E_R(R/\p))$
  is the zero map, for all $i\geq 0$. Thus, $\Tor^R_i(R/I,E_R(R/\p))=0$ for all $i\geq 0$. Since for each $i\geq 0$, the torsion functor
  $\Tor^R_i(R/I,-)$ commutes with the direct sums it follows that $\Tor^R_i(R/I,E_R(M))=0$, for all $i\geq 0$. So, from the exact
   sequence $(3.2.1)$ it follows that $\Tor^R_i(R/I,E_R(M)/M)=0$, for all $i\geq 0$. Hence, by applying the inductive hypothesis to the
    $R$-module $E_R(M)/M$ we have $$ \Ext^{i+1}_R(R/I,M)\simeq\Ext^{i}_R(R/I,E_R(M)/M)=0\,\,\,\,{\rm for}\,\,0\leq i \leq n-1.$$ So, we have $\Ext^n_{R}(R/I,M)=0$. This completes the inductive step.\qed\\

(ii)$\Rightarrow$(i) Assume the opposite.  Then there is an integer $j\geq 0$ such that $\Tor^R_{j}(R/I,M)\neq 0$.
Let $\p \in \Supp \Tor^R_{j}(R/I,M)$. Then localizing at $\p$, without loss of generality, we may assume that $(R,\m)$
 is a Noetherian local ring. Let $T:=D(M)$, where $D(-)$
denotes the Matlis dual functor. Then by the adjointness we have $$ \Tor^R_i(R/I,T)\simeq D(\Ext^i_R(R/I,M))=0,$$
for all $i\geq 0$. Therefore, by the previous part of the proof we have  $\Ext^i_R(R/I,T)=0,$ for each $i\geq 0$.
So, by the adjointness we have $$ D(\Tor^R_j(R/I,M))\simeq\Ext^j_R(R/I,T)=0,$$which implies that $\Tor^R_j(R/I,M)=0$.
This is the desired contradiction.\qed\\

\begin{lem}
\label{3.3}
Let $R$ be a Noetherian ring and $I$ be an ideal of $R$. Let $E$ be an injective $R$-module and $K$ be a submodule of $E$. Then $$\frac{L}{\Gamma_I(L)}\otimes_R \frac{R}{I}=0,\,\,\,\,\,\,\,{\rm where}\,\, L:=\frac{E}{K}.$$
\end{lem}
\proof  In view of \cite[Proposition 2.1.4]{BS} the $R$-module $E_1:=\Gamma_I(E)$ is injective.
Therefore, there exists an injective submodule $E_2$ of $E$ such that $E_1+E_2=E$ and $E_1\cap E_2=0$.
 Since $\frac{E_1+K}{K}$ is a submodule of $\Gamma_I(L)$ it follows that the $R$-module $\frac{L}{\Gamma_I(L)}$
 is a homomorphic image of the $R$-module $E_2$. So, in order to prove the assertion it is
 enough to prove that $E_2\otimes_R \frac{R}{I}=0$. So, we must prove that $E_2=I E_2$. Since $\Gamma_I(R)E_2\subseteq \Gamma_I(E_2)=0$, it follows that $\Gamma_I(R)E_2=0$. Therefore, we have $$E_2=(0:_{E_2}\Gamma_I(R))\simeq \Hom_R(R/\Gamma_I(R),E_2).$$
 On the other hand, by \cite[Lemma 2.1.1]{BS} there is an exact sequence $$0\longrightarrow R/\Gamma_I(R)
  \stackrel{a} \longrightarrow R/\Gamma_I(R),$$for some element $a\in I$, which effecting the $R$-linear
  exact functor $\Hom_R(-,E_2)$ induces the exact sequence $$ \Hom_R(R/\Gamma_I(R),E_2)\stackrel{a} \longrightarrow \Hom_R(R/\Gamma_I(R),E_2) \longrightarrow0.$$ Therefore, we have $\Hom_R(R/\Gamma_I(R),E_2)=a \Hom_R(R/\Gamma_I(R),E_2)$. Hence, we have
\begin{eqnarray*}
E_2/aE_2&\simeq&E_2\otimes_RR/Ra\\&\simeq&\Hom_R(R/\Gamma_I(R),E_2)\otimes_RR/Ra\\&\simeq&\Hom_R(R/\Gamma_I(R),E_2)/a\Hom_R(R/\Gamma_I(R),E_2)\\&\simeq&0.
\end{eqnarray*}
So, we have $E_2=aE_2$ and hence $E_2=IE_2$.\qed\\

\begin{lem}
\label{3.4}
Let $R$ be a Noetherian ring and let $I\subseteq J$ be two ideals of $R$.
Let $M$ be an $R$-module and $t\geq 2 $ be an integer such that $\Tor^R_j(R/J,H^i_I(M))=0$ for all $i>t$ and all $j\geq0$. Then we have $\Tor^R_j(R/J,H^t_I(M))=0$ for $j=0,1$.
\end{lem}
\proof Let $$0\longrightarrow M \stackrel{\varepsilon}\longrightarrow E_0
\stackrel{f_0}\longrightarrow E_1 \stackrel{f_1}\longrightarrow E_2 \stackrel{f_2}\longrightarrow\cdots $$
be a minimal injective resolution for $M$. Set $K_i:={\rm ker} f_i$ for $i\geq0$.
By splitting this minimal injective resolution to some short exact sequences we get the isomorphisms
$$H^t_I(M)\simeq H^{t-1}_I(K_1)\simeq H^{t-2}_I(K_2)\simeq\cdots\simeq H^2_I(K_{t-2})
\simeq H^1_I(K_{t-1})\simeq H^1_I(K_{t-1}/\Gamma_I(K_{t-1})).$$
Set $N:=K_{t-1}/\Gamma_I(K_{t-1})$. Then we have $\Gamma_I(N)=0$ and $$K_{t-1}=
{\rm ker} f_{t-1}={\rm im} f_{t-2} \simeq E_{t-2}/{\rm ker} f_{t-2}=E_{t-2}/K_{t-2}.$$
So, there exists a submodule $K$ of $E:=E_{t-2}$ such that $N\simeq E/K$. By \cite[Remark 2.2.7]{BS}
 there is an exact sequence $$0  \longrightarrow N  \longrightarrow D_I(N)  \longrightarrow H^1_I(N)
   \longrightarrow0.\,\,\,\,\,(3.4.1)$$ For each $i\geq 2$ we have $$H^i_I(N)\simeq H^i_I(K_{t-1})\simeq H^{t+i-1}_I(M).$$
Therefore, by the hypothesis we have $$\Tor^R_j(R/J,H^i_I(N))\simeq \Tor^R_j(R/J,H^{t+i-1}_I(M))=0$$
for all $i\geq 2$ and all $j\geq0$.
 Thus, by Proposition 3.2 we have $\Ext^j_R(R/J,H^i_I(N))=0$ for all $i\geq 2$ and all $j\geq0$.
  The exact sequence $(3.4.1)$ yields the isomorphisms $$ H^i_I(D_I(N))\simeq H^i_I(N),$$ for $i\geq 2$.
  So, we have $$\Ext^j_R(R/J,H^i_I(D_I(N)))\simeq \Ext^j_R(R/J,H^i_I(N))=0$$ for all $i\geq 2$ and all $j\geq0$.
   Also, by \cite[Corollary 2.2.8]{BS} we have $H^i_I(D_I(N))=0$, for $i=0,1$. So $\Ext^j_R(R/J,H^i_I(D_I(N)))=0$
    for all $i\geq 0$ and all $j\geq0$. Therefore, by \cite[Lemma 2.1]{GBA} we have $\Ext^j_R(R/J,D_I(N))=0$,
    for all integers $j\geq 0$. Thus by Proposition 3.2 we have $\Tor^R_{i}(R/J,D_I(N))=0$, for all integers $i\geq 0$.
     On the other hand, in view of Lemma 3.3 we have $N\otimes_R R/I=0$ and hence it follows from the hypothesis
     $I \subseteq J$ that $N\otimes_R R/J=0$. Hence, using the long exact sequence induced by the exact sequence
     $(3.4.1)$ it follows that $\Tor^R_{j}(R/J,H^1_I(N))=0$, for $j=0,1$. Therefore, $\Tor^R_{j}(R/J,H^t_I(M))\simeq\Tor^R_{j}(R/J,H^1_I(N))=0$, for $j=0,1$.\qed\\

\begin{cor}
\label{3.5}
Let $R$ be a Noetherian ring and $I\subseteq J$ be two ideals of $R$. Let $K$ be an $R$-module such that $2\leq{\rm cd}(I,K)=t$. Then $\Tor^R_i(R/J,H^t_I(K))=0$, for $i=0,1$.
\end{cor}
\proof The assertion follows from Lemma 3.4. \qed\\

\begin{prop}
\label{3.6}
 Let $R$ be a Noetherian ring, $I$ be an
ideal of $R$ and $M$ be an $R$-module. Then the following conditions
are equivalent:
\begin{itemize}
 \item[(i)] $\Tor^R_{i}(R/I,M)=0$, for all $i\geq 0$.
 \item[(ii)] $\Tor^R_{i}(R/I,M)=0$, for all $0\leq i\leq \cd(I,R)$.
\end{itemize}
\end{prop}
\proof (i)$\Rightarrow$(ii) Is clear.\\

(ii)$\Rightarrow$(i) Assume the opposite. Then there is an integer $j>\cd(I,R)$
such that $\Tor^R_{j}(R/I,M)\neq 0$. Let $\p \in \Supp \Tor^R_{j}(R/I,M)$. Then
localizing at $\p$, without loss of generality, we may assume that $(R,\m)$ is a
 Noetherian local ring. Let $T:=D(M)$, where $D(-)$
denotes the Matlis dual functor. Then by the adjointness we have $$\Ext^i_R(R/I,T)\simeq D(\Tor^R_i(R/I,M))=0,$$
for all $0\leq i\leq \cd(I,R)$. Now, by \cite[Theorem 2.9]{AM} we have $H^i_I(T)=0$, for all $0 \leq i \leq \cd(I,R)$.
 Hence, $H^i_I(T)=0$, for all integers $i\geq 0$. So, in view of \cite[Theorem 2.9]{AM} we have $\Ext^i_R(R/I,T)=0,$
 for all integers $i\geq 0$. Consequently, by the adjointness, we have $$ D(\Tor^R_i(R/I,M))\simeq\Ext^i_R(R/I,T)=0,$$
 for all $i\geq 0$. Thus, $\Tor^R_i(R/I,M)=0$ for all $i\geq 0$, which is a contradiction.\qed\\

\begin{thm}
\label{3.7}
  Let $R$ be a Noetherian ring and $I$ be an
ideal of $R$. Then the following conditions
are equivalent:
\begin{itemize}
 \item[(i)]  ${\rm cd}(I,R)\leq 1$.
 \item[(ii)] The ideal transform functor $D_I(-)$ is exact.
 \item[(iii)] For every $R$-module $M$, if $\Tor^R_{i}(R/I,M)=0$ for $i=0,1$, then $\Tor^R_{i}(R/I,M)=0$, for all integers $i\geq 0$.
 \item[(iv)] For every $R$-module $M$, if $\Ext^i_{R}(R/I,M)=0$ for $i=0,1$, then $\Ext^i_{R}(R/I,M)=0$, for all integers $i\geq 0$.
\end{itemize}
\end{thm}
\proof (i)$\Leftrightarrow$(ii) See \cite[Lemma 6.3.1]{BS}.\\

(i)$\Rightarrow$(iii) Follows from Proposition 3.6.\\

(iii)$\Rightarrow$(i) Assume the opposite. Then we have
$\cd(I,R)\geq 2$. Let $t=\cd(I,R)$. Then by Corollary 3.5
for the $R$-module $M:=H^t_I(R)$ we have $\Tor^R_{i}(R/I,M)=0$
for $i=0,1$. So, by the hypothesis we have $\Tor^R_{i}(R/I,M)=0$, for all integers
$i\geq 0$. Therefore, by Proposition 3.2 we have $\Ext^i_{R}(R/I,M)=0$, for all
 integers $i\geq 0$. Therefore $(0:_{H^t_I(R)}I)\simeq\Hom_R(R/I,H^t_I(R))=0$.
 Hence $H^t_I(R)=0$, which is a contradiction.\\

(i)$\Rightarrow$(iv) Let $M$ be an $R$-module such that $\Ext^i_{R}(R/I,M)=0$
for $i=0,1$. Then, by \cite[Theorem 2.9]{AM} we have $H^i_I(M)=0$ for $i=0,1$.
On the other hand by \cite[Lemma 6.3.1]{BS} we have $H^i_I(M)=0$, for all integers $i\geq 2$.
Therefore, $H^i_I(M)=0$, for all integers $i\geq0$.  Hence, it follows from \cite[Theorem 2.9]{AM}
 that $\Ext^i_R(R/I,M)=0$ for all integers $i\geq 0$.\\

(iv)$\Rightarrow$(i) By \cite[Corollary 2.2.8]{BS} we have $H^i_I(D_I(R))=0$, for $i=0,1$.
  Hence, $\Ext^i_R(R/I,D_I(R))=0$ for $i=0,1$ by \cite[Theorem 2.9]{AM}.
  So, by the hypothesis we have $\Ext^i_R(R/I,D_I(R))=0$, for all integers $i\geq 0$.
  Hence, by Proposition 3.2 it follows that $R/I\otimes_R D_I(R)=0$ and so $D_I(R)=ID_I(R)$.
  Now the assertion follows from \cite[Lemma 6.3.1 and Proposition 6.3.5]{BS}.\qed\\

\begin{thm}
\label{3.8}
 Let $R$ be a Noetherian ring and let $I=Ra_1+Ra_2$ be an
ideal of $R$. Then the following conditions
are equivalent:
\begin{itemize}
 \item[(i)]  ${\rm cd}(I,R)=2$.
 \item[(ii)] There exists an $R$-module $M$, such that $\Tor^R_{i}(R/I,M)=0$ for $i=0,1$ and $\Tor^R_{2}(R/I,M)\neq0$.
 \item[(iii)] There exists an $R$-module $M$, such that $\Ext^i_{R}(R/I,M)=0$ for $i=0,1$ and $\Ext^2_{R}(R/I,M)\neq0$.
\end{itemize}
\end{thm}
\proof
(i)$\Rightarrow$(ii) If ${\rm cd}(I,R)=2$, then for the $R$-module
$M:=H^2_I(R)$ by Corollary 3.5 we have $\Tor^R_{i}(R/I,M)=0$ for $i=0,1$.
 Now, we claim that $\Tor^R_{2}(R/I,M)\neq0$. Assume the opposite. Then, it follows
 from Proposition 3.6 that $\Tor^R_{i}(R/I,M)=0$, for all integers $i\geq 0$. Then it
  follows from Proposition 3.2 that $\Ext^i_{R}(R/I,M)=0$, for all integers $i\geq 0$.
   Therefore, $\Hom_R(R/I,H^2_I(R))=0$, which implies that $H^2_I(R)=0$. This is a contradiction.\\

(ii)$\Rightarrow$(i) Under the given hypothesis it follows from Theorem 3.7
 that ${\rm cd}(I,R)\geq2$. On the other hand, by \cite[Theorem 3.3.1]{BS} we
 have ${\rm cd}(I,R)\leq2$. So, we have ${\rm cd}(I,R)=2$.\\

(i)$\Rightarrow$(iii) If ${\rm cd}(I,R)=2$, then for the $R$-module $M:=D_I(R)$,
by \cite[Corollary 2.2.8]{BS} we have $H^i_I(M)=0$, for $i=0,1$.  Hence, by \cite[Theorem 2.9]{AM}
 we have $\Ext^i_R(R/I,M)=0$ for $i=0,1$.
Moreover, by \cite[Remark 2.2.7]{BS} there is an exact sequence $$0  \longrightarrow R/\Gamma_I(R) \longrightarrow D_I(R)  \longrightarrow H^1_I(R)  \longrightarrow0,$$ which induces the isomorphisms $$H^2_I(M)\simeq H^2_I(R/\Gamma_I(R)) \simeq H^2_I(R)\neq 0.$$ Now, if $\Ext^2_R(R/I,M)=0$, then  $\Ext^i_R(R/I,M)=0$, for $i=0, 1, 2,$. So, it follows from \cite[Theorem 2.9]{AM} that $H^i_I(M)=0$ for $i=0,1,2$,
which is a contradiction, because $H^2_I(M)\neq 0$. So, for the $R$-module $M:=D_I(R)$ we have $\Ext^i_{R}(R/I,M)=0$ for $i=0,1$ and $\Ext^2_{R}(R/I,M)\neq0$.\\

(iii)$\Rightarrow$(i) Under the given hypothesis it follows from Theorem 3.7
 that ${\rm cd}(I,R)\geq2$. On the other hand, by \cite[Theorem 3.3.1]{BS}
 we have ${\rm cd}(I,R)\leq2$. So, we have ${\rm cd}(I,R)=2$.\qed\\

\begin{rem} \label{3.9} For a given proper ideal $I$  of a Noetherian ring $R$,
 there are some other well known equivalent conditions for $\cd (I, R)\leq1$.
 For instance, see \cite[Lemma 6.3.1 and Proposition 6.3.5]{BS} and for more
 properties of such ideals see \cite{B}.
\end{rem}

\section{Vanishing of the Bass and Betti numbers of local cohomology modules}

In this section we prepare some vanishing conditions for some of the Bass and Betti numbers of special local cohomology modules. Also, we give a formula for the cohomological dimension of special ideals over Noetherian domains.

\begin{thm}
\label{4.1}
  Let $R$ be a Noetherian ring and let $I\subseteq J$ be two ideals of $R$ such that $\cd(J,R)=1$. Then the following statements hold:
\begin{itemize}
 \item[(i)]  For every $R$-module $M$ we have $\Ext^j_R(R/J,H^i_I(M))=0$, for all integers $i\geq 2$ and $j\geq 0$.
 \item[(ii)] For every $R$-module $M$ and every finitely generated $R$-module $N$ with $\Supp N \subseteq V(J)$ we have $\Ext^j_R(N,H^i_I(M))=0$, for all integers $i\geq 2$ and $j\geq 0$.
 \item[(iii)] For every $R$-module $M$ and each prime ideal $\p\in V(J)$ we have $\mu^j(\p,H^i_I(M))=0$, for all integers $i\geq 2$ and $j\geq 0$. $($Here $\mu^j(\p,H^i_I(M))$ denotes the j-th Bass number of the $R$-module $H^i_I(M)$ with respect to $\p$$)$.
 \item[(iv)]  For every $R$-module $M$ and each prime ideal $\p\in V(J)$ we have $\beta_j(\p,H^i_I(M))=0$, for all integers $i\geq 2$ and $j\geq 0$. $($Here $\beta_j(\p,H^i_I(M))$ denotes the j-th Betti number of the $R$-module $H^i_I(M)$ with respect to $\p$$)$.
\end{itemize}
\end{thm}
\proof (i) Let $M$ be an arbitrary $R$-module. In order to prove the assertion we may assume $\cd(I,M)=t\geq 2$. By Corollary 3.5 we have $\Tor^R_j(R/J,H^t_I(M))=0$ for $j=0,1$.  Now, it follows from Proposition 3.6 that $\Tor^R_i(R/J,H^t_I(M))=0$, for all integers $i\geq 0$. Therefore, by Proposition 3.2 we have $\Ext^i_{R}(R/J,H^t_I(M))=0$, for all integers $i\geq 0$.  Now, if $t\geq 3$ then by Lemma 3.4 we have $\Tor^R_j(R/J,H^{t-1}_I(M))=0$ for $j=0,1$.  So, it follows from Proposition 3.6 that $\Tor^R_i(R/J,H^{t-1}_I(M))=0$, for all integers $i\geq 0$. Therefore, by Proposition 3.2 we have $\Ext^i_{R}(R/J,H^{t-1}_I(M))=0$, for all integers $i\geq 0$. Proceeding in the same way we see that $\Ext^j_R(R/J,H^i_I(M))=0$, for all integers $i\geq 2$ and $j\geq 0$.\\

(ii) Using \cite[Lemma 2.2]{AB1} follows from (i).\\

(iii) Follows from (ii).\\

(iv) Follows from (ii) using Proposition 3.2. \qed\\

\begin{thm}
 \label{4.2}
 Let $R$ be a Noetherian domain and let $I$ and $J$ be two non-zero proper
 ideals of $R$ such that $\cd(J,R)=1$ and $J\not\subseteq \Rad (I)$. Then we have $$\cd(I\cap J,R)=\max\{i\in \Bbb{Z}\,\,:\,\,\Supp H^i_I(R)\not\subseteq V(J)\}.$$
\end{thm}
\proof Since by the hypothesis we have $J\not\subseteq \Rad (I)$ it follows that  $$J\not\subseteq \bigcap_{\p\in\mAss_R R/I}\p,$$which implies that $J\not\subseteq \q$ for some $\q\in \mAss_R R/I$. Assume that $k:=\height \q$. Then by Grothendieck's
 Non-vanishing Theorem we have $$0 \neq H^k_{\q R_{\q}}(R_{\q})=H^k_{I R_{\q}}(R_{\q})\simeq (H^k_I(R))_{\q},$$ which implies that $\q \in \Supp H^k_I(R)$. Hence we have $\Supp H^k_I(R)\not\subseteq V(J)$ and so we have $$k\in\{i\in \Bbb{Z}\,\,:\,\,\Supp H^i_I(R)\not\subseteq V(J)\}.$$
In particular, we have $\{i\in \Bbb{Z}\,\,:\,\,\Supp H^i_I(R)\not\subseteq V(J)\}\neq\emptyset.$\\

Now, set $\ell:= \max\{i\in \Bbb{Z}\,\,:\,\,\Supp H^i_I(R)\not\subseteq V(J)\}$ and $t:=\cd(I\cap J,R)$.
 Then as by \cite[Corollary 3.3.3]{BS} we have $\ell\leq{\rm ara}(I)$ and $t\leq{\rm ara}(I\cap J)$ it
 follows that $0\leq\ell<\infty$ and $0\leq t<\infty$. By the Mayer-Vietoris exact sequence for
  each integer $i>t$ we have the exact sequence $$H^{i}_{I+J}(R)\longrightarrow H^{i}_I(R)\oplus H^i_J(R) \longrightarrow H^i_{I\cap J}(R),$$
  which gives
 the exact sequence $$H^{i}_{I+J}(R)\longrightarrow H^{i}_I(R)\oplus H^i_J(R) \longrightarrow0$$and hence
 we have $$\Supp H^i_I(R)\subseteq \Supp H^{i}_{I+J}(R) \subseteq V(I+J)\subseteq V(J).$$ So, it is clear
 that $\ell \leq t$. On the other hand, since by the hypothesis $R$ is a domain and $I\cap J\neq 0$, it
 follows that $\ell \geq 1$ and $t \geq 1$. Hence, if $t=1$, then it is clear that $\ell=1=t$. Now,
 assume that $t\geq2$ and $\ell<t$.  Then, we have $\Supp H^t_I(R)\subseteq V(J)$. Also, by the Mayer-Vietoris exact sequence
    we have the exact sequence $$H^{t}_I(R)\oplus H^t_J(R) \longrightarrow H^t_{I\cap J}(R)\longrightarrow H^{t+1}_{I+J}(R),$$
    which implies that $$\Supp H^t_{I\cap J}(R)\subseteq [\Supp H^t_J(R) \cup \Supp H^t_I(R) \cup \Supp H^{t+1}_{I+J}(R)]\subseteq V(J).$$
    So, the non-zero $R$-module $H^t_{I\cap J}(R)$ is $J$-torsion and hence we have $$\Hom_R(R/J ,H^t_{I\cap J}(R))\neq 0.$$
    But, by Theorem 4.1 we have $$\Hom_R(R/J ,H^t_{I\cap J}(R))=0,$$ which is a contradiction. So, we have $\ell=t$, whenever $t\geq 2$. \qed\\

\begin{prop}
 \label{4.3}
  Let $(R,\m)$ be a Noetherian local ring and let $I$ and $J$ be two proper ideals of $R$
  such that $\cd(J,R)=1$. Let $k\geq 2$ be an integer such that $\Supp H^k_I(R)\not\subseteq V(J)$.
  Then, the $R$-module $H^k_{I\cap J}(R)$ is not $I\cap J$-cofinite. In particular, $H^k_{I\cap J}(R)\neq 0$.
\end{prop}
\proof
  By the Mayer-Vietoris exact sequence we have the exact sequence
  $$H^k_{I+J}(R)\longrightarrow H^{k}_I(R)\oplus H^k_J(R) \longrightarrow H^k_{I\cap J}(R),$$
  which considering that fact that $\Supp H^k_{I+J}(R) \subseteq V(I+J)\subseteq V(J)$
  implies that $\Supp H^k_{I\cap J}(R)\not\subseteq V(J).$ (Note that by the hypothesis
  we have $\Supp H^k_I(R)\not\subseteq V(J)$.) In particular, we have $H^k_{I\cap J}(R)\neq 0$.
  In order to prove the assertion, assume the opposite and assume that $\dim \Supp H^k_{I\cap J}(R)=d$.
   Then, in view of \cite[Theorem 2.9]{M} we have $H^d_{\m}(H^k_{I\cap J}(R))\neq 0$. On the other hand,
    by Theorem 4.1 we have $\Ext^j_R(R/\m,H^k_{I\cap J}(R) )=0$, for all integers $j\geq 0$.
    Hence it follows from \cite[Theorem 2.9]{AM} that $H^j_{\m}(H^k_{I\cap J}(R))=0$, for all
    integers $j\geq 0$. Therefore, we have $H^d_{\m}(H^k_{I\cap J}(R))=0$, which is a contradiction. \qed\\

\begin{rem}
\label{4.4}
  There are examples of Noetherian local rings $(R,\m)$ with proper ideals $I$,
  for which $\cd (I, R)=1$ and $\ara(I)\geq 2$. For instance, the following example is given by Hellus and St\"uckrad in \cite{HS1}.
\end{rem}

\begin{ex}
\label{4.5}
Let $k$ be a field and let $S=k[[x,y,z,w]]$,  where $x, y, z, w$ are independent
indeterminacies  over $k$. Let $f=xw-yz$, $g=y^3-x^2z$ and $h=z^3-w^2y$. Let $R=S/fS$
and $I=(f,g,h)S/fS$. Then $R$ is a Noetherian local ring of dimension $3$ with maximal
 ideal $\m=(x,y,z,w)S/fS$. Also, for the ideal $I$ of $R$, we have $\cd(I,R)=1$ and
  $\ara(I)= 2$. $($See \cite[Remark 2.1(ii)]{HS1}$)$.
\end{ex}

\begin{center}
{\bf Acknowledgments}
\end{center}

The authors are deeply grateful to the referee for a very careful
reading of the manuscript and many valuable suggestions.
Also, we would like to thank Professors Reza Naghipour and Kamran
Divaani-Aazar for their careful reading of the first draft and many
helpful suggestions.


\end{document}